\documentclass[12pt,a4paper]{article} \textheight=23cm
\usepackage{makeidx}
\usepackage{amssymb}
\usepackage{graphics,graphicx}
\usepackage[T1]{fontenc}
\usepackage{color}

\def\EE{\mathbb{E}}
\def\ZZ{\mathbb{Z}}
\def\RR{\mathbb{R}}
\def\TT{\mathbb{T}}
\def\TT{\mathbb{T}}
\def\NN{\mathbb{N}}
\def\PP{\mathbb{P}}

\def\exp{{\rm exp}}

\def \mod {{\rm \ mod \ }}

\def\Proof {\vskip -2mm {\underline {Proof}}}
\def\fdm {\hfill\break \boxit{2pt} {}}
\def\build#1_#2^#3{\mathrel{\mathop{\kern 0pt#1}\limits_{#2}^{#3}}}

\font\fivegoth=eufm5 \font\sevengoth=eufm7 \font\tengoth=eufm10
\newfam\gothfam
\scriptscriptfont\gothfam=\fivegoth \textfont\gothfam=\tengoth
\scriptfont\gothfam=\sevengoth

\def\boxit#1#2{\setbox1=\hbox{\kern#1{#2}\kern#1}%
\dimen1=\ht1 \advance\dimen1 by #1 \dimen2=\dp1 \advance\dimen2 by#1
\setbox1=\hbox{\vrule height\dimen1 depth\dimen2\box1\vrule}%
\setbox1=\vbox{\hrule\box1\hrule}%
\advance\dimen1 by .4pt \ht1=\dimen1 \advance\dimen2 by .4pt
\dp1=\dimen2 \box1\relax}
\def \fdm {\boxit{3pt} {}}

\newtheorem{hypo}{Hypothesis}[section]
\newtheorem{ppty}[hypo]{Property}
\newtheorem{prop}[hypo]{Proposition}
\newtheorem{thm}[hypo]{Theorem}
\newtheorem{lem}[hypo]{Lemma}
\newtheorem{defi}[hypo]{Definition}

\newtheorem{notas}[hypo]{Notations}

\begin{document}
\parindent=0mm

\markboth{Conze \& Le Borgne}{Limit law for some modified ergodic sums}

\title{Limit law for some modified ergodic sums}

\date{}

\maketitle

\begin{centerline}{Jean-Pierre CONZE, St\'ephane LE BORGNE}\end{centerline}
\vskip 4mm
\begin{centerline}{Irmar, Universit\'e de Rennes 1\footnote{Irmar, UMR CNRS
6625, Universit\'e de Rennes 1, Campus de Beaulieu, 35042 Rennes Cedex, France.}}\end{centerline}

\vskip2cm

\begin{abstract}
An example due to Erd\H{o}s and Fortet shows that, for a lacunary
sequence of integers $(q_n)$ and a trigonometric polynomial
$\varphi$, the asymptotic distribution of ${1\over \sqrt n}
\sum_{k=0}^{n-1} \varphi(q_k x)$ can be a mixture of gaussian laws.
Here we give a generalization of their example interpreted as the
limiting behavior of some modified ergodic sums in the framework of
dynamical systems.
\end{abstract}



\vskip 8mm
\section*{Introduction}

\vskip 3mm Let $(q_n)$ be a lacunary sequence of natural integers.
The stochastic-like behavior of the sums $\sum_{k=0}^{n-1}
\varphi(q_k x)$ for $\varphi$ a regular 1-periodic real function has
been the subject of several works (Fortet, Kac, Salem, Zygmund,...).
Before the second war and in the forties, different particular cases
were treated for which the Central Limit Theorem (CLT) can be shown.
The fact that the CLT is not always satisfied in its standard form
was already noticed by Fortet et Erd\H{o}s. Gaposhkin \cite{Ga66},
Berkes \cite{Be76}, recently Berkes and Aisleitner \cite{AiBe08}
gave arithmetic conditions on the sequence $(q_n)$ which imply the
CLT.

The counter-example of Fortet and Erd\H{o}s is very simple. Let us take
$q_n = 2^n -1$ and $\varphi(x) = \cos(2\pi x) +\cos(4\pi x)$. Then
the limit law of the distribution of the normalized sums
$n^{-{1\over 2}} \sum_{k=0}^{n-1} \varphi((2^k-1) x)$ is not the
gaussian one, but a mixture of gaussian laws, explicited in
\cite{Ka47}, \cite{Ka49}\footnote{The proof announced by Kac, as
well as an article by Erdos, Ferrand, Fortet and Kac on the sums
$\sum_{k=0}^{n-1} \varphi(q_n x)$ mentioned in \cite{Ka49}, did not
appear as far as we know. For a proof based on a result of Salem and
Zygmund, \cite{SaZy48}) for lacunary sequences see Berkes and
Aisleitner \cite{AiBe08}. Here we give a slightly different proof
and generalizations of it.}.

This fact is related to the arithmetic properties of the sequence
$2^n -1$. But it can also be interpreted from other points of view.
It can be viewed as a consequence of non ergodicity for stationary
martingales which gives asymptotically a mixture of gaussian laws.

As we will show in Section \ref{groupAct}, one can equally view it
as a special case of a general phenomenon for isometric perturbation
of dynamical systems of hyperbolic type.

For instance, let $A\in SL(d,\ZZ)$ be a matrix without eigenvalue
root of the unity, let $B$ be a matrix with integral coefficients
and $\det B\neq 0$. Denoting by $\lambda$ the Lebesgue measure on
the $d$-dimensional torus $\TT^d$, we have, for every centered
{H\"older} function $\varphi$ on $\TT^d$, for every $t \in \RR$,
\begin{eqnarray*} \lambda\{x: \ \frac{1}{\sqrt{n}}
\sum_0^{n-1}\varphi((A^k-B)x) < t \} \longrightarrow \int_{\TT^d} \,
\frac{1}{\sqrt{2\pi}\sigma_y} \int_{-\infty}^t
\exp(\frac{-s^2}{2\sigma_y^2})ds\ d\lambda(y),
\end{eqnarray*}
where $\sigma_y$ is the asymptotic variance of the translated
function: $\varphi_y(x) := \varphi(x+y)$.

\vskip 3mm In Section \ref{multiMix} we will also use a different
method, based on a property of multiple decorrelation, to extend the
results to a large class of chaotic dynamical systems and their
modified ergodic sums.

\vskip 3mm In what follows $(X,d)$ will be a metric space with its
Borel $\sigma$-algebra ${\cal B}$ and a probability measure $\mu$ on
${\cal B}$, $T $ a measure preserving transformation on $(X, {\cal
B}, \mu)$, and $\varphi$ a real function on $X$ with some
regularity. We denote by $S_n \varphi$ (or $S_n (\varphi)$) the
ergodic sums of $\varphi$:
$$S_n \varphi(x) := \sum_{k=0}^{n-1} \varphi( T ^k x).$$

\section{Generalization of an example of Fortet and Erd\H{o}s}

\vskip 3mm In order to explain the method on a simple example, we
begin with the counter-example of Fortet and Erd\H{o}s mentioned in the
introduction. Then we indicate how it can be extended to modified
ergodic sums for expanding maps of the interval.

\subsection{The example of Fortet and Erd\H{o}s}

The space $X$ is the circle $\RR/\ZZ$, $\mu$ is the Lebesgue measure
and $T$ is the transformation $x \rightarrow 2x {\rm \ mod \ 1 }$.

\begin{notas} We denote by $R_n \varphi(x,y)$ the translated modified ergodic sums
\begin{eqnarray}
R_n \varphi(x,y) := \sum_{k=0}^{n-1} \varphi(2^k x- x + y).
\end{eqnarray}
\end{notas}

Let $\varphi$ be an {H\"older}ian function on $\RR/\ZZ$ such that
$\int_0^1 \varphi \, d\mu = 0$. We denote by $c_p(\varphi)$ its
Fourier coefficient of order $p \in \ZZ$. We have:
$$n^{-1} \|S_n\varphi\|_2^2 =  \sum_{p\not = 0} [|c_p(\varphi)|^2
+ 2 \sum_{1\leq j \leq n-1} (1-{j\over n}) \, c_{2^jp}(\varphi)
\overline c_p(\varphi)].$$ The variance of $\varphi$ is well defined
and given by
$$ \sigma^2(\varphi) := \lim_n {\|S_n\varphi\|_2^2 \over n} =
\sum_{p\not = 0} [|c_p(\varphi)|^2
+ 2 \sum_{j \geq 1} \, c_{2^jp}(\varphi) \overline c_p(\varphi)].$$

\vskip 3mm For $\varphi(x) = \cos(2\pi x) + \cos(4\pi x)$, we have
the following convergence for every $t \in \RR$:
\begin{eqnarray} \lim_n \mu\{x :  {1\over \sqrt{n}} \sum_{k=0}^{n-1} \varphi((2^k-1) x)
\leq t \} = {1\over \sqrt{2\pi}} \int_0^1 (\int_{-\infty}^{t/ |\cos
y|} e^{-s^2/2} ds)\  dy,\label{tclFEr}
\end{eqnarray}
or in terms of characteristic function: \begin{eqnarray}\lim_n
\EE(e^{it {1\over \sqrt{n}} \sum_{k=1}^{n} \varphi((2^k -1) .)}) =
\int_0^1 e^{-{1\over 2}(\cos y)^2 \,t^2} \ dy.
\label{caracx}\end{eqnarray}

Before we give generalizations of this result, for the reader's
convenience we recall a proof of (\ref{caracx}) (cf. \cite{AiBe08})
based on the following general statement:

\begin{lem} \label{phi} Let $(Z_n)$ be a sequence of real random variables
on $[0,1]$. Let $\cal{L}$ be a probability distribution on $\RR$,
with characteristic function $\Phi(t) = \int e^{itx} \ {\cal
L}(dx)$. The following conditions are equivalent:

a) for every probability density $\rho$, the sequence $(Z_n)$ under
the measure $\rho \mu$ converges in distribution to $\cal{L}$;
\hfill \break b) for every interval $I \subset [0,1]$,
\begin{eqnarray}\lim_n {1\over \mu(I)} \mu\{x \in I: Z_n(x) \leq t\}
= {\cal L}(]-\infty, t]), \ \forall t \in \RR;
\label{int-dens}\end{eqnarray} \hfill \break c) for every Riemann
integrable function $\psi$, the sequence $(\psi Z_n)$ converges in
distribution to a limit distribution with characteristic function
$\int_0^1 \Phi(\psi(y)\,t) \ dy$.

In particular if ${\cal L} = {\cal N}(0,1)$, under the previous
conditions the sequence $(\varphi Z_n)$ converges in distribution to
a limit distribution whose characteristic function is $\int_0^1
e^{-{1\over 2}\varphi^2(y)\, t^2} \, dy$.
\end{lem}
\Proof \ \ \ Assume b). Let $\psi$ be a step function, $\psi =
\sum_{j=0}^p c_j 1_{[a_j, \,a_{j+1}[}$, with $a_0 = 0 < a_1 < ... <
a_{p+1} =1$. We have:
$$\EE_\mu(e^{it\psi(.) Z_n(.)}) = \sum_j \int_{a_j}^{a_{j+1}} e^{it c_j
Z_n(.)} \ d\mu \rightarrow \sum_j \mu(I_j) \Phi(c_jt) = \int_0^1
\Phi(\psi(y) \,t) \ dy.$$ The general case c) follows by approaching
$\psi$ by step functions.

\vskip 3mm Conversely, let $\rho = \sum_{j=0}^p c_j 1_{[a_j,
\,a_{j+1}[}$, with $\rho \geq 0$ and $\int_0^1 \rho \ dx = 1$. Under
Condition b) or c), we have: $$\EE_{\rho \mu}(e^{it Z_n(.)}) =
\sum_j c_j \int_{a_j}^{a_{j+1}} e^{it Z_n(.)} \ d\mu \rightarrow
[\sum_j c_j (a_{j+1} - a_j)] \Phi(t) = \Phi(t).$$ As above we obtain
the general case by approximation. Condition a) follows. \fdm

\vskip 3mm Salem and Zygmund proved in \cite{SaZy48} the CLT with
Condition b) of Lemma \ref{phi} for $\varphi =\cos$ and any lacunary
sequence $(q_k)$. The convergence (\ref{caracx}) follows from their
result, from Lemma \ref{phi}, and from the trigonometric identity
($n \geq 2$):
\begin{eqnarray*}&&\sum_{k=1}^{n}[\cos(2\pi (2^k -1)\, x) + \cos(4\pi (2^k -1)\, x)] \\
&&= \cos(2\pi x) + \cos(2\pi (2^{n+1} -2)\, x) + 2 \cos(\pi x)
\sum_{k=2}^{n} \cos(2\pi (2^k  - 3/2) \, x).
\end{eqnarray*}

\vskip 2mm Now we give an analogous result for more general
functions. The method of proof is slightly different from the
previous one and will be applied to various examples in the sequel
of the paper.

\vskip 3mm First of all we need a well known improved version of the
CLT for regular functions. In the special case of this section, it
can be proved using the classical method of quasi-compact operator.
We first sketch the idea for the transformation $x \rightarrow 2x
{\rm \ mod \ 1 }$ on $\RR/\ZZ$. For this map, the dual
(Perron-Frobenius) operator $P$ is given by
$$P\varphi(x) = {1\over 2} [\varphi({x\over 2}) + \varphi({x\over 2} + {1\over 2})].$$

Let $\rho$ be a probability density on $X = \RR/\ZZ$. Then for every
integer $\ell \geq 0$, we have:
\begin{eqnarray*}
&&|\int_X e^{i{t \over  \sqrt{n}} S_n \varphi (x)} \ \rho(x) \ dx -
\int_X e^{i{t \over  \sqrt{n}} S_n \varphi (x)} \ \ dx| \\ &=&
|\int_X [e^{i{t \over  \sqrt{n}} S_n \varphi (x)} - e^{i{t \over
\sqrt{n}} S_n \varphi (T ^\ell x)}] \ \rho(x) \ dx + \int_X e^{i{t
\over \sqrt{n}}
S_n \varphi (x)} \ (P^\ell \rho(x) - 1) \ dx|\\
&\leq& 2 {|t|\over \sqrt{n}} \ell \|\varphi\|_\infty + \|P^\ell \rho
- 1\|_1.
\end{eqnarray*}

If $\rho$ is a function of bounded variation or a {H\"older} function,
then $\|P^\ell \rho - 1\|_1$ converges to 0 when $\ell$ tends to
$\infty$ with an exponential rate. If $\rho$ is only in $L^1$,
convergence holds a priori without rate.

\vskip 3mm Therefore, in the CLT, we can replace the Lebesgue
measure by a measure which is absolutely continuous with respect to
the Lebesgue measure when the density is regular. It allows us to
apply Lemma \ref{phi} or similar results. Actually this principle
holds for dynamical systems in a very general situation as we will
see later. Applied here to the density $\mu(D)^{-1} 1_{D-y} \mu $
where $D$ is an interval in $[0,1[$, it yields the following result:

\begin{lem} \label{sigmay} For every $y$, every interval $D$,
every regular 1-periodic function $\varphi$,
$$\mu(D)^{-1}\mu\{x: n^{-{1\over 2}} \sum_{k=0}^{n-1} \varphi(2^k x +y) \leq t
{\rm \ and \ } x \in D-y \} \rightarrow {1\over \sqrt{2\pi}}
\int_{-\infty}^{t/ \sigma(\varphi_y)} e^{-s^2/2} \, ds,$$ where
$\sigma(\varphi_y)$ is the  asymptotic variance of the translated
function: $\varphi_y(x) := \varphi(x+y)$.
\end{lem}

We will also use a property of regularity of the variance with
respect to translations:
\begin{lem} If $C(\varphi):=\sum_p |p| \, |c_p(\varphi) | < +\infty$, then we have:
\begin{eqnarray} n^{-{1\over 2}} \|R_n \varphi(.,y) - R_n \varphi(.,y')
\|_2 \leq C(\varphi) |y-y'|.\label{reg-y}
\end{eqnarray}
\end{lem} \Proof \ \ For every $p \not = 0$, we have ${1\over n} \int
|\sum_{k=0}^{n-1} e^{2\pi i 2^k px}|^2 \ dx =1.$

\vskip 3mm Writing $\varphi(2^k x -x +y) = \sum_{p \not = 0}
c_p(\varphi) e^{2i\pi p (y-x)} e^{2i\pi p 2^k x}$, the following
uniform bound holds:
\begin{eqnarray*} {1\over \sqrt{n}}\|R_n \varphi(.,y) - R_n \varphi(.,y')\|_2
&\leq& |y-y'|\sum_p |p|\,|c_p(\varphi) | \|{1\over \sqrt{n}} \|S_n
e^{2\pi i p.} \|_2 \\
&=& |y-y'|\sum_p |p| \, |c_p(\varphi) | = C(\varphi) |y-y'|.
\end{eqnarray*} \fdm

\vskip 3mm By convention in the sequel, when the variance $\sigma$
is $0$, $\int_{-\infty}^{t/\sigma} e^{-u^2/2} \ du$ is interpreted
as the repartition function of the limit law, the Dirac mass at 0.

\begin{thm} If $\varphi$ satisfies $\sum_p |p| \, |c_p(\varphi) | < +\infty$, we
have:
\begin{eqnarray}\mu\{x :  {1\over \sqrt{n}} \sum_{k=0}^{n-1} \varphi(2^k x -x) \leq t \}
\rightarrow {1\over \sqrt{2\pi}} \int_0^1 (\int_{-\infty}^{t/
\sigma(\varphi_y)} e^{-s^2/2} \, ds)\  dy.\label{tclx}
\end{eqnarray}
\end{thm}
\Proof \ \  It is enough to check the convergence when $t$ is a
continuity point of the limit distribution. By integrating with
respect to $y$ and applying Lebesgue theorem, it follows from Lemma
\ref{sigmay}:
\begin{eqnarray*}
& &{1\over \mu(D)} \mu\{(x,y) : n^{-{1\over 2}} \sum_{k=0}^{n-1} \varphi(2^k x +y) \leq t
{\rm \ and \ } x \in D-y \} \\
& &\ \ \ \ \ \ \ \ \ \ \ \rightarrow {1\over \sqrt{2\pi}}
\int_0^1 (\int_{-\infty}^{t/ \sigma(\varphi_y)} e^{-s^2/2} \, ds)\
dy.
\end{eqnarray*}
The change of variable $(x,y) \rightarrow (x, y+x)$ leaves the
measure $dx \times dy$ invariant. We get that the difference:
\begin{eqnarray} & &{1\over \mu(D)}\mu\{(x,y) : n^{-{1\over 2}} \sum_{k=0}^{n-1} \varphi(2^k
x -x +y) \leq t {\rm \ and \ } y \in D \} \\
& &\ \ \ \ \ \ \ \ \ \ \ - {1\over \sqrt{2\pi}}
\int_0^1 (\int_{-\infty}^{t/ \sigma(\varphi_y)} e^{-s^2/2} \, ds)\
dy \label{intervD}
\end{eqnarray}
converges to 0. For $\gamma >0$, we obtain by (\ref{reg-y}):
\begin{eqnarray*}
&&\mu\{x : {1\over \sqrt{n}} R_n \varphi(x, 0) \leq t \} \\ &&\leq
\mu\{x : {1\over \sqrt{n}} R_n \varphi(x,y) \leq t+\gamma \}+ \mu\{x
: {1\over \sqrt{n}} |R_n \varphi(x, 0) -R_n \varphi(x, y)| \geq
\gamma \},
\end{eqnarray*}

hence:
\begin{eqnarray}
&&\mu\{x : {1\over \sqrt{n}} R_n \varphi(x, 0) \leq t \} \nonumber \\
&&\leq {1\over \mu(D)} \int_X \mu\{x : {1\over \sqrt{n}} R_n
\varphi(x,y) \leq t+\gamma \} \ 1_D(y) \ dy + C(\varphi)^2 {\delta^2
\over \gamma^2}, \label{maj1} \\
&&{1\over \mu(D)} \int_X \mu\{x : {1\over \sqrt{n}} R_n \varphi(x,y)
\leq t-\gamma \} \ 1_D(y) \ dy - C(\varphi)^2 {\delta^2 \over
\gamma^2} \nonumber \\ &&\leq \mu\{x : {1\over \sqrt{n}} R_n
\varphi(x, 0) \leq t \}.\label{min1}
\end{eqnarray}

Let $\varepsilon >0$. First we take $\gamma$ such that
$$|{1\over \sqrt{2\pi}} \int_0^1 (\int_{-\infty}^{t \pm \gamma / \sigma(\varphi_y)}
e^{-s^2/2} \, ds)\ dy - {1\over \sqrt{2\pi}} \int_0^1
(\int_{-\infty}^{t/ \sigma(\varphi_y)} e^{-s^2/2} \, ds)\ dy | <
\varepsilon,$$ then $\delta$ such that $C(\varphi)^2 {\delta^2 \over
\gamma^2} < \varepsilon$ and finally $n_0$ such that the difference
in (\ref{intervD}) (with $t\pm\gamma$ in place of $t$) is less then $\varepsilon$.

\vskip 3mm Applying (\ref{maj1}) and (\ref{min1}), we get
$$|\mu\{x :  {1\over \sqrt{n}} \sum_{k=0}^{n-1} \varphi(2^k x -x) \leq t \}
- {1\over \sqrt{2\pi}} \int_0^1 (\int_{-\infty}^{t/
\sigma(\varphi_y)} e^{-s^2/2} \, ds)\  dy | < 6 \varepsilon.$$ \fdm

\vskip 3mm \goodbreak {\it Non-degeneracy of the limit distribution}

For a regular function $\varphi$, the set $\{y: \sigma(\varphi_y) =
0 \}$ is closed, since $\sigma(\varphi_y)$ is continuous as a
function of $y$. This set coincides with $E(\varphi,T ):= \{y:
\varphi_y {\rm \ is \ a \ } T {\rm -coboundary}\}$. Since $\varphi$
is regular, if $\varphi_y = \psi_y\circ T - \psi_y$ for a function
$\psi_y$, then  $\psi_y$ is also regular.

\vskip 3mm If the set $E(\varphi,T )$ has measure 1, it coincides
with $[0,1]$. The functions $\varphi_y$ vanish for every $y$ on the
fixed point of $T$, which implies that $\varphi$ is identically
zero. Therefore the limiting distribution in (\ref{tclx}) is non
degenerated when $\varphi$ is non identically zero.

\vskip 3mm
\goodbreak
\subsection{First generalizations, expanding maps}

\vskip 4mm  {\bf Expanding maps}

The previous example can be extended to other dynamical systems or
"sequential dynamical systems" in different directions. For instance
we can consider a lacunary sequence of positive integers $(q_n)$ and
a modified version of the corresponding ergodic sums, like
$\sum_{q=0}^{n-1} \varphi(q_k x -x)$.

We can also take the class of expanding maps on the interval, for example
the $\beta$-transformations, for which the spectral properties of
the transfer operator can be used.

\vskip 3mm We will not develop these specific extensions, but we
will consider the following general framework. The example of Erd\H{o}s
and Fortet can be view as a special case of the following general
construction.

\vskip 3mm Let us consider a dynamical system $(X, T , \mu)$, a
space ${\cal F}$ of real valued functions on $X$, and a map $\theta
: x \rightarrow \theta_x$ from $X$ into the set of Borel maps from
$X$ to $X$ which preserve ${\cal F}$ under composition. Thus, for
each $x \in X$, a map $\theta_x$ (also denoted by $\theta(x)$) is
given.

\vskip 3mm Suppose that for $\varphi \in {\cal F}$  the sums
$\sum_{k=0}^{n-1} \varphi(T ^k x)$ after normalization have a
distribution limit, for instance convergence in distribution of
$n^{-{1\over 2}} \sum_{k=0}^{n-1} \varphi(T ^k x)$ toward a gaussian
law. Now let us consider the modified sums:
\begin{equation} \sum_{k=0}^{n-1} \varphi(\theta_x(T ^k x)). \label{gen-form}
\end{equation}

In the previous section, the system was $x \rightarrow 2 x \mod 1$
on the circle, ${\cal F}$ the space of {H\"older} functions or the space
of functions $\varphi$ satisfying $\sum_p |p| \, |c_p(\varphi) | <
+\infty$, and the map $\theta_x : y \rightarrow y-x$. In the sequel
of the paper we will describe different cases where a limit law for
the modified sums (\ref{gen-form}) can be obtained.

A simple situation is the following. Let us consider a map $\theta$
taking a finite number of values. More precisely, suppose that there
is a finite partition $(A_j, j \in J)$ of $X$ in measurable sets
$A_j$ such that $\theta_x  = \theta_j$ on $A_j$, where $\theta_j$ is
a {H\"older}ian map from $X$ to $X$.

Applying a result of Zweim\"uller \cite{Zw07}, we have:
\begin{eqnarray*}
\EE [e^{it {1 \over \sqrt n} \sum_{k=0}^{n-1} \varphi(\theta(.)(T ^k
.))}] &=& \sum_{j\in J} \mu(A_j)\, \int \, e^{it {1 \over \sqrt n}
\sum_{k=0}^{n-1} \varphi(\theta_j(T ^k
x))} \ \mu(A_j)^{-1} 1_{A_j} \, d\mu \\
&&\rightarrow \sum_{j\in J} \mu(A_j)\, e^{-{1 \over
2}\sigma^2(\varphi\,\circ \,\theta_j) \ t^2}.
\end{eqnarray*}

We would like to extend such a convergence to maps $\theta$ taking a
continuum of values like in the example of Erd\H{o}s and Fortet. In the
next section, this will be done when there is a compact group $G$
acting on $X$ and when $\theta$ is a regular map with values in $G$.
If we denote the variance $\sigma^2(\varphi\,\circ \,\theta_x)$ by
$\sigma_x^2$, our aim is to prove that the limit distribution of the
normalized modified sums (\ref{gen-form}) has for characteristic
function:
$$\int_X \, e^{-{1 \over 2}\sigma_x^2 \, t^2} \ d\mu (x).$$

\vskip 3mm \section{Generalization to group actions
\label{groupAct}}

\subsection{A general result}

Let $(X,T,\mu)$ be a dynamical system. Suppose that a compact group
$K$ acts on $X$ and preserves the measure $\mu$. Denote by $m$ the
Haar measure on $K$, and by $(k,x) \rightarrow kx$ the action of $K$
on $X$.

Let $\theta \ : \ X \rightarrow K$ be a Borel function. The modified
sums that we consider here have the form
\begin{eqnarray*} \sum_{k=0}^{n-1} \varphi(\theta(x) \, T ^k x).
\end{eqnarray*}

Let $\mu_\theta$ the measure on $K$ defined by $\mu_\theta(B) =
\mu(\theta^{-1} B)$, for $B$ Borel set in $K$. For a function
$\varphi$ in $L^2(\mu)$, we introduce the two following properties:

\vskip 3mm \begin{ppty} \label {RR1} (Central limit theorem with
density for $\varphi(k.)$) There exists $\sigma_k \geq 0$ such that,
for every $\rho$ density with respect to $\mu$ on $X$, for every
interval $A$ in $\RR$, every $k\in K$,
\begin{eqnarray} (\rho \mu) \{x \, : \, \frac{1}{\sqrt{n}} S_n \varphi(k x)
\in A  \} \longrightarrow \frac{1}{\sqrt{2\pi} \sigma_k} \int_A
\exp(- {1\over 2} \frac{s^2}{\sigma_k^2})\ ds. \label{sigmak}
\end{eqnarray}\end{ppty}

If $\sigma_k=0$ in (\ref{sigmak}), the limit law is the
Dirac measure at 0.

\begin{ppty} \label {RR2} (Continuity of the variance of the modified
sums with respect to the translations)
$$\|\frac{1}{\sqrt{n}} \sum_{\ell =0}^{n-1}\varphi(\theta(\cdot)^{-1}k \; T^\ell \cdot) -
\frac{1}{\sqrt{n}} \sum_{\ell =0}^{n-1}\varphi(\theta(\cdot)^{-1}
T^\ell \cdot)\|_2\leq C\sqrt{n}d(k,e).$$ \end{ppty}

\begin{prop} \label{thetak} If the function $\varphi$ satisfies \ref{RR1}
and \ref{RR2}, the following convergence holds:
\begin{eqnarray*} \mu\{x: \ \frac{1}{\sqrt{n}}
\sum_{\ell=0}^{n-1}\varphi(\theta(x)\,T^\ell x) \in A \}
\longrightarrow \int_{K}\frac{1}{\sqrt{2\pi} \sigma_k} \int_A \exp(-
{1\over 2} \frac{s^2}{\sigma_k^2})ds\ d\mu_{\theta}(k).
\end{eqnarray*} \end{prop}

\Proof \ \  First, let us fix $D$ a compact neighborhood of the
identity in $K$. For every element $k$ of $K$, the regularity of the
action of $K$ and Property \ref{RR1} for the function
$\varphi_k(\cdot)=\varphi(k\cdot)$ give the convergence
\begin{eqnarray*}
& &\mu\{x \, : \,  \frac{1}{\sqrt{n}} S_n\varphi_k(x) <u, \
\theta(x)^{-1}k\in D  \} \\
& &\ \ \ \ \ \ \ \ \ \longrightarrow \frac{1}{\sqrt{2\pi}
\sigma_k} \int_{-\infty}^u \exp(- {1\over 2}
\frac{s^2}{\sigma_k^2})ds\ \mu\{x\ : \ \theta(x)^{-1}k\in D \},\end{eqnarray*}
where $\sigma_k^2$ is the asymptotic variance associated to the
function $\varphi_k$.

By taking the integral over $K$, we obtain:
\begin{eqnarray*}
{ } & &
(\mu\otimes m) \{(x,k)\ : \ \frac{1}{\sqrt{n}} S_n\varphi_k(x) <u,\  \theta(x)^{-1}k\in D  \}\\
& &\ \ \ \ \ \longrightarrow \int_{K}\frac{1}{\sqrt{2\pi} \sigma_k}
\int_{-\infty}^u \exp(- {1\over 2} \frac{s^2}{\sigma_k^2})ds\
\mu(\{x \ : \ \theta(x)^{-1}k\in D\})\ dk.
\end{eqnarray*}
The measure $\mu\otimes m$ is preserved by the change of variable
 $x=x,\ k'=\theta(x)^{-1}k$. So, by dividing by
$m(D)$, we get:
\begin{eqnarray*}
{ } & &{{1}\over{m(D)}}\mu\otimes m\{(x,k')\ : \ \frac{1}{\sqrt{n}}
\sum_{\ell=0}^{n-1}\varphi(\theta(x)\, k' \, T^\ell x) <u ,\  k'\in D  \}\\
& &\ \ \ \ \longrightarrow \int_{K}\frac{1}{\sqrt{2\pi} \sigma_k}
\int_{-\infty}^u \exp(- {1\over 2} \frac{s^2}{\sigma_k^2})ds\
{{\mu(\{x : \ \theta(x)^{-1}k\in D\})}\over{m(D)}}\ dk.
\end{eqnarray*}
The measure
$${{\mu(\{x : \ \theta(x)^{-1}k\in D\})}\over{m(D)}}\ dk$$
tends to $\mu_\theta$ (the image of $\mu$ by $\theta$) when the
diameter of $D$ tends to 0. Property \ref{RR2} now allows us to
conclude. Let $\varepsilon$ be a positive number. Let
$R_n\varphi(x,k)$ be the sum
$$R_n\varphi(x,k)= \sum_{\ell=0}^{n-1}\varphi(\theta(x)kT^\ell x).$$
We have
\begin{eqnarray*}
& &\mu(\{x\ : \ {{1}\over{\sqrt{n}}}R_n\varphi(x,e)\leq u\})\\
& &\ \ \ \ \ \leq \mu(\{x\ : \ {{1}\over{\sqrt{n}}}R_n\varphi(x,k)
\leq u+\varepsilon\})\\
& &\ \ \ \ \ \ \ \ \ \ \ \ \ +\mu(\{x\ : \ {{1}\over{\sqrt{n}}}|R_n \varphi(x,k)-R_n\varphi(x,e)|>\varepsilon\})\\
& &\ \ \ \ \ \leq \mu(\{x\ : \ {{1}\over{\sqrt{n}}}R_n\varphi(x,k)
\leq u+\varepsilon\})+{{Var(R_n\varphi(\cdot,k)-R_n\varphi(\cdot,e))}\over{n\varepsilon^2}}\\
& &\ \ \ \ \ \leq \mu(\{x\ : \ {{1}\over{\sqrt{n}}}R_n\varphi(x,k)
\leq u+\varepsilon\})+{{C^2nd(k,e)^2}\over{n\varepsilon^2}}.
\end{eqnarray*}

Denoting by $\delta(D)$ the diameter of $D$, the average taken on
$D$ yields:
\begin{eqnarray*}
& &\mu(\{x\ : \ {{1}\over{\sqrt{n}}}R_n\varphi(x,e)\leq u\})\\
& &\ \ \ \ \ \leq {{1}\over{m(D)}}\mu\otimes m\{(x,k): \
\frac{1}{\sqrt{n}} \sum_{l=0}^{n-1}\varphi(\theta(x)kT^lx)
<u+\varepsilon ,\  k'\in D \}+
{{C^2\delta(D)^2}\over{\varepsilon^2}}.
\end{eqnarray*}
and, because of the above convergence,
\begin{eqnarray*}
& &\limsup_{n\rightarrow\infty}\mu(\{x\ : \ {{1}\over{\sqrt{n}}}R_n\varphi(x,e)\leq u\})\\
& &\ \ \ \ \ \leq \int_{K}\frac{1}{\sqrt{2\pi} \sigma_k}
\int_{-\infty}^{u+\varepsilon} \exp(\frac{-s^2}{2\sigma_k^2})ds\
{{\mu(\{x \ : \ \theta(x)^{-1}k\in D\})}\over{m(D)}}\
dk+{{C^2\delta(D)^2}\over{\varepsilon^2}},\\
\end{eqnarray*}
thus, letting $\delta(D)$ tends to 0,
\begin{eqnarray*}
\limsup_{n\rightarrow\infty}\mu(\{x\ : \
{{1}\over{\sqrt{n}}}R_n\varphi(x,e)\leq u\}) \ \leq \
\int_{K}\frac{1}{\sqrt{2\pi} \sigma_k}
\int_{-\infty}^{u+\varepsilon} \exp(\frac{-s^2}{2\sigma_k^2})ds\
d\mu_{\theta}(k).
\end{eqnarray*}
Similarly we have
\begin{eqnarray*}
\liminf_{n\rightarrow\infty}\mu(\{x\ : \ {{1}\over{\sqrt{n}}}
R_n\varphi(x,e)\leq u\}) \ \geq \int_{K}\frac{1}{\sqrt{2\pi}
\sigma_k} \int_{-\infty}^{u-\varepsilon} \exp(-{1\over 2} \frac{s^2}
{\sigma_k^2})ds\ d\mu_{\theta}(k).
\end{eqnarray*}
Therefore, excepted maybe for $u = 0$ if $\sigma_k= 0$ for a set of
positive $\mu_\theta$-measure of elements $k$, we have the
convergence:
\begin{eqnarray*} \lim_{n\rightarrow\infty}\mu(\{x\ : \
{{1}\over{\sqrt{n}}}R_n\varphi(x,e)\leq u\}) \
=\int_{K}\frac{1}{\sqrt{2\pi} \sigma_k} \int_{-\infty}^{u}
\exp(-{1\over 2} \frac{s^2}{\sigma_k^2})ds\ d\mu_{\theta}(k).
\end{eqnarray*}
\fdm

\vskip 3mm A typical situation in which one might apply this result
is when the central limit theorem holds for regular functions and
the action of $K$ is regular. Let $(X,d)$ be a metric space. For a
real number $\eta>0$, on the space of {H\"older} continuous functions of
order $\eta$ we define the $\eta$-variation and the $\eta$-{H\"older}
norm by
\begin{eqnarray}
[\varphi]_{\eta}=\sup_{x \neq y}
{{|\varphi(x)-\varphi(y)|}\over{d(x,y)^\eta}}, \
\|\varphi\|_{\eta}=\|\varphi\|_\infty+[\varphi]_{\eta}.\label{etaNorm}
\end{eqnarray}

We say that the action of $K$ is {H\"older} continuous if, for every
$k\in K$, $x\rightarrow kx$ is {H\"older} continuous. For many chaotic
systems it has been proved that the central limit theorem holds for
{H\"older} continuous functions. If the action of $K$ is {H\"older}
continuous, then the central limit theorem holds for $\varphi(k.)$.
Moreover, because of the theorem of Eagleson \cite{Ea76,Zw07}, the theorem is true for measures absolutely continuous
with respect to $\mu$. In this case a {H\"older} continuous function
$\varphi$ satisfies Property \ref{RR1}.

\vskip 3mm \begin{defi} {\rm We say that $(X,T,\mu)$ is {\it
summably mixing} if, for every $\eta>0$, there exist $C>0$ and a
summable sequence $(\alpha_n)$ such that, for every centered
$\eta$-{H\"older} continuous functions $\varphi$, $\psi$, one has :
$$
|\int_X\varphi \circ T^n\ \psi\ d\mu|\leq
 C\|\varphi\|_\eta
 \|\psi\|_2\alpha_n.
$$
}\end{defi} If $(X,T,\mu)$ is summably mixing and $\varphi$ is a
centered $\eta$-{H\"older} continuous function then there exists $C>0$
such that
\begin{eqnarray}
\|\sum_{\ell=0}^{n-1}\varphi(T^\ell\cdot)\|_2\leq
C\|\varphi\|_\eta^{1/2}\|\varphi\|_2^{1/2}\, \sqrt{n}.\label{l2}
\end{eqnarray}
\vskip 3mm
\begin{prop} \label{groupes} Let $(X,T,\mu)$ be a dynamical system where $X$
is a Riemannian manifold.  Assume that there is a measure preserving
action of a compact Lie group $K$ on $(X, \mu)$ which is $C^\infty$.
Assume that the central limit theorem holds for differentiable
functions on $X$ and that $(X,T,\mu)$ is summably mixing. Then for
every $C^\infty$ function $\varphi$ on $X$, for every $t\in\RR$,
\begin{eqnarray*}
\mu\{x: \ \frac{1}{\sqrt{n}} \sum_0^{n-1}\varphi(\theta(x)T^kx) <t
\} \longrightarrow \int_{K}\frac{1}{\sqrt{2\pi} \sigma_k}
\int_{-\infty}^t \exp(-{1\over 2} \frac{s2}{\sigma_k2}) \, ds\
d\mu_{\theta}(k).
\end{eqnarray*}
\end{prop}

\Proof \ \ We will use Fourier analysis on $K$. We briefly summarize
what we need. For more details see \cite{Bo82}.

The action of $K$ on $X$ defines a unitary representation $U$ of $K$
on $L^2(\mu)$ by $k\longmapsto \varphi(k^{-1}x)$ which can be
decomposed as a sum of irreducible representations. Let $\hat{K}$ be
the set of the equivalence classes of irreducible representations of
$K$ and $\delta$ be an element of $\hat{K}$.

Let us fix a base  $R$ of the root system of ${\cal K}$ the Lie
algebra of $K$. Let us call $W$ the associated Weyl chamber. To each
irreducible representation of $K$ is uniquely associated a linear
form belonging to a lattice in $W$: the dominant weight of the
representation. Let $\delta$ be an element of $\hat{K}$ and let
$\gamma_{\delta}$ be the corresponding dominant weight.

\vskip 3mm The Weyl formula gives the dimension $d_{\delta}$ of the
irreducible representation associated to $\delta$ as a function of
$\gamma$:
$$d_{\delta}=  \prod_{\alpha\in R_+}\frac{\langle
\alpha,\gamma_{\delta}+\rho\rangle} {\langle \alpha,\rho\rangle},$$
where $R_+$ is the set of positive roots and $\rho$ the half sum of
the positive roots. \vskip 3mm For every $\delta\in \hat{K}$, let
$\xi_{\delta}$ be the character of $\delta$,
$\chi_{\delta}=d_{\delta}\xi_{\delta}$, and
\begin{eqnarray}
P_{\delta}=U(\overline{\chi_{\delta}})=d_{\delta}\int_K
\overline{\xi_{\delta}(k)} \, U(k)\ dk.\label{proj}
\end{eqnarray}
The operator $P_{\delta}$ is the projection of $L^2(\mu)$ on the
isotypic part ${\cal F}_{\delta} := P_{\delta}(L^2(\mu))$. We have
the decomposition
$$L^2(\mu)=\bigoplus_{\delta\in\hat{K}}{\cal F}_{\delta}.$$
For a given vector $v$ in $L^2(\mu)$ let $v_\delta := P_\delta v$.
An element $v$ of ${\cal F}_{\delta}$ is $K$-finite:
\begin{eqnarray}
{\rm dim\ Vect}Kv\leq d_{\delta}^2. \label{Kfini}
\end{eqnarray}
One says that $v$ is $C^\infty$ if the map $k\mapsto U(k)v$ is
$C^\infty$ . One defines the derived representation of $U$ on the
space of $C^\infty$  elements; it is a representation of the Lie
algebra ${\cal K}$ of $K$ and that one can be extended to a
representation of the universal enveloping algebra of ${\cal K}$. We
use the same later $U$ to denote these three representations.

\vskip 3mm Let $X_1,\ldots,X_{n}$ be an orthonormal basis for an
invariant scalar product on ${\cal K}$. The operator
$\Omega=1-\sum_{i=1}^nX_i^2$ belongs to the center of the universal
enveloping algebra of ${\cal K}$. So, by Schur's lemma, if
$\mu_{\delta}$ is a representation of the type $\delta$, there
exists $c_{\delta}$ such that
$\mu_{\delta}(\Omega)=c_{\delta}\mu_{\delta}(1)$.

The operators $\Omega(X_i)$ being hermitian, $c_{\delta}$ is
positive. One can show ({\sl cf.} \cite{Bo82}) that there exists a
scalar product $Q$ such that $c_{\delta} =
Q(\gamma_{\delta}+\rho)-Q(\rho)$.

If $v$ is $C^\infty$, one has
$$P_{\delta}U(\Omega)v=c_{\delta}P_\delta v=c_{\delta}v_\delta,$$
thus, for every non negative integer $m$, for every $\delta$ in
$\hat{K}$, one has
$$v_\delta= c_{\delta}^{-m}\ (U(\Omega^m)v)_\delta,$$
with large $c_{\delta}$ for large $\gamma_\delta$. From this
equality and the definition (\ref{proj}) of $P_\delta$, one deduces
that
$$\|v_\delta\|_\infty\leq
{{d_\delta^2}\over{c_\delta^m}}\|U(\Omega^m)v\|_\infty.$$

In particular, the series $\sum_{\delta\in\hat{K}}v_\delta$
converges uniformly to $v$. Therefore we can write, for every $x\in
X$,
\begin{eqnarray*}
\varphi(k^{-1}x)&=&\sum_{\delta\in\hat{K}}U(k)(\varphi_\delta)(x), \\
\varphi(\theta(x)kx)&=& \sum_{\delta\in\hat{K}}
U((\theta(x)k)^{-1})(\varphi_\delta)(x).
\end{eqnarray*}

We want to study the quantity:
$$\|\frac{1}{\sqrt{n}} \sum_{\ell=0}^{n-1}\varphi(\theta(\cdot)kT^\ell \cdot)-
 \frac{1}{\sqrt{n}} \sum_{\ell =0}^{n-1}\varphi(\theta(\cdot)T^\ell \cdot)\|_2.$$
With the notations introduced above we can write:
\begin{eqnarray*}
\sum_{\ell =0}^{n-1}\varphi(\theta_xkT^\ell x)-
\frac{1}{\sqrt{n}} \sum_{\ell =0}^{n-1}\varphi(\theta_xT^\ell x)
&=&\sum_{\ell =0}^{n-1}((U(k^{-1})-Id)U(\theta_x^{-1})\varphi)(T^\ell x)\\
&=&\sum_{\delta\in\hat{K}}\sum_{\ell
=0}^{n-1}((U(k^{-1})-Id)U(\theta_x^{-1})\varphi_\delta)(T^\ell x).
\end{eqnarray*}
Since $v_\delta$ is $K$-finite, there exists a finite set of
$C^\infty$  functions $\{\varphi_{\delta,j}, \, \ j=1\ldots q\}$
(with $\varphi_{\delta,j}=U(k_j) \varphi_{\delta}$ for some $k_j$
and $q\leq d_\delta^2$ because of (\ref{Kfini})) and uniformly
bounded functions $u_{\delta,j}$ such that
$$
U(\theta(x)^{-1})\varphi_\delta=\sum_{j=1}^qu_{\delta,j}(x)\varphi_{\delta,j}.
$$
Furthermore one has
\begin{eqnarray}
\|(U(k^{-1})-Id)v_\delta\|\leq
C\|\gamma_\delta\|^rd(k,e)\|v_\delta\|,\label{moder}
\end{eqnarray}
where $\|\gamma_\delta\|$ denotes a usual norm of a point in a
lattice (Inequality (\ref{moder}) corresponds to the fact that the
norm of the operator $(U(k^{-1})-Id)$ on ${\cal F}_\delta$ is a
"moderately increasing" function of $\delta$ (\cite{Bo82}, p.
82-83)).

Thus there exists $C>0$, and regular functions $w_{\delta,i}$ with
$|w_{\delta,i}(k)|\leq C\|\gamma_\delta\|^rd(k,e)$ such that
$$(U(k^{-1})-Id)U(\theta(x)^{-1})\varphi_\delta
=\sum_{i,j=1}^qw_{\delta,i}(k)u_{\delta,j}(x)\varphi_{\delta,j},$$

The triangular inequality now gives
\begin{eqnarray*}
& &|\sum_{l=0}^{n-1}((U(k^{-1})-Id)U(\theta(x)^{-1})\varphi_\delta)(T^lx)|\\
& &\ \ \ \ \ \ \ \leq \sum_{i,j=1}^q|w_{\delta,i}(k)u_{\delta,j}(x)||\sum_{l=0}^{n-1}\varphi_{\delta,i}(T^lx)|.\\
\end{eqnarray*}
But the norms $\|\varphi_{\delta,i}\|$ are equal to
$\|\varphi_{\delta}\|$ and the norms $\|\varphi_{\delta,i}\|_\eta$
are bounded by $d_\delta^2\|\varphi\|_\eta$ (observe that with
formula (\ref{proj}) we can control the regularity of
$\varphi_\delta=P_\delta(\varphi)$). Thus, thanks to the boundness
of the variance (\ref{l2}) applied to $\varphi_{\delta,j}$, we
obtain
\begin{eqnarray*}
& &\|\sum_{l=0}^{n-1}((U(k^{-1})-Id)U(\theta(x)^{-1})\varphi_\delta)(T^lx)\|_2\\
& &\leq
C\sum_{i,j=1}^q\|\gamma_\delta\|^{r}d(k,e)\|\sum_{l=0}^{n-1}\varphi_{\delta,j}(T^lx)\|_2
\leq C\sum_{i,j=1}^q\|\gamma_\delta\|^{r}d(k,e)\|\varphi_{\delta,j}\|_\eta^{1/2}\|\varphi_{\delta,j}\|_2^{1/2}\, \sqrt{n}\\
& & \displaystyle{ \leq
Cq^2\|\gamma_\delta\|^{r}d(k,e)\|\varphi_\delta\|_2^{1/2}d_\delta\|\varphi\|_\eta^{1/2}\sqrt{n}
\leq
Cd^5_\delta\|\gamma_\delta\|^{r}\|\varphi_\delta\|_2^{1/2}\|\varphi\|_\eta^{1/2}d(k,e)\sqrt{n}.}
\end{eqnarray*}
We then again use the triangular inequality to obtain
\begin{eqnarray*}
& &\|\sum_{\delta\in\hat{K}}\sum_{l=0}^{n-1}((U(k^{-1})-Id)U(\theta(\cdot)^{-1})\varphi_\delta)(T^l\cdot)\|_2\\
& &\ \ \ \ \ \ \ \leq
C(\sum_{\delta\in\hat{K}}d_\delta^5\|\gamma_\delta\|^{r}\|\varphi_\delta\|_2^{1/2})
\|\varphi\|_\eta^{1/2} \ \sqrt{n}\ d(k,e).
\end{eqnarray*}
Since $\varphi$ is $C^\infty$ the series
$\sum_{\delta\in\hat{K}}d_\delta^5\|\gamma_\delta\|^{r}\|\varphi_\delta\|_2^{1/2}$
converges and we get the inequality
$$\|\frac{1}{\sqrt{n}} \sum_{l=0}^{n-1}\varphi(\theta(\cdot)kT^l\cdot)-
\frac{1}{\sqrt{n}}
\sum_{l=0}^{n-1}\varphi(\theta(\cdot)T^l\cdot)\|_2\leq
C(\varphi)d(k,e).$$ \fdm

\vskip 3mm
\subsection{Examples}

\vskip 3mm {\bf Automorphisms of the torus} ({\it cf.} \cite{Le60,LB99})

Let $A\in SL(d,\ZZ)$ be a matrix without eigenvalue root of the
unity. It defines an ergodic automorphism of the $d$-dimensional
torus $\TT^d$ for the Lebesgue measure $\lambda$. Let $B$ be a
matrix with integral coefficients. For every centered $C^\infty$
function $\varphi$, for every $t \in \RR$, we have:
\begin{eqnarray*} \lambda\{x: \ \frac{1}{\sqrt{n}}
\sum_0^{n-1}\varphi((A^k-B)x) < t \} \longrightarrow
\int_{\TT^d}\frac{1}{\sqrt{2\pi} \sigma_y} \int_{-\infty}^t
\exp(\frac{-s^2}{2\sigma_y^2})ds\ d\lambda_{B}(y).
\end{eqnarray*}
If $\det B\neq 0$ then $\lambda_{B}=\lambda$.

\vskip 3mm {\it Non-degeneracy of the limit}

\vskip 2mm If $T $ is a hyperbolic automorphism of the torus
$\TT^d$, then for any {H\"older}ian $\varphi$, the set $E(\varphi, T )
:= \{y: \varphi(.+y) {\rm \ is \ a \ } T {\rm -coboundary}\}$ is
closed because $\varphi(.+y)$ is a coboundary if and only
$\sigma(\varphi_y)=0$ and $y\mapsto \sigma(\varphi_y)$ is a
continuous function (this is a consequence of the mixing properties
of $T$ ; see section 4.2 below). On the other hand if $\varphi(.+y)
{\rm \ is \ a \ } T {\rm -coboundary}$ then it is a coboundary
inside the set of {H\"older} continuous functions: there exists a {H\"older}
continuous function $\psi_y$ such that $\varphi(.+y)=T
\psi_y-\psi_y$. Thus if $B$ is surjective and $\sigma_y=0$
$\lambda$-almost surely, then $\sigma_y=0$ everywhere and
$\varphi(0+y)=T \psi_y(0)-\psi_y(0)=\psi_y(0)-\psi_y(0)=0$. So, if
$B$ is surjective, unless $\varphi=0$, the limit law is not
degenerated.

\vskip 3mm If $T$ is an ergodic  non-hyperbolic automorphism of the
torus $\TT^d$ then one has to reinforce the regularity hypothesis on
$\varphi$ in order to apply the second part of the previous
reasoning: if $\varphi$ is $d$-times differentiable, has a {H\"older}
continuous $d^{th}$-differential and is a measurable coboundary then
it is a coboundary in the space of {H\"older} continuous functions
\cite{Ve86,LB99}).

\vskip 3mm {\bf Automorphisms on nilmanifolds} ({\it cf.}
\cite{CoLB02})

Let $X$ be the 3-dimensional nilmanifold defined as the homogeneous
space $N/\Gamma$, where $N$ is the Heisenberg group of triangular
matrices $$\pmatrix {
1 & x_1 & z \cr 0 & 1 & x_2 \cr 0 & 0 & 1}$$ and $\Gamma$ the discrete subgroup of
integral points in $N$. Let $\mu$ be the $N$-invariant measure on
$N/\Gamma$ induced by the Haar measure on $N$. We identify $N$ and
$\RR^3$ equipped with the law
$$(x_1,x_2,z) . (x_1',x_2',z') = (x_1+x_1',x_2+x_2', z+z'+x_1x_2' - x_1'x_2).$$

Let $$A =\pmatrix {
a & b \cr c & d \cr },$$ be a hyperbolic matrix in $Sl(2, \ZZ)$. It defines a
transformation $T$ on $N/\Gamma$ by
$$T : (x_1,x_2,z)\Gamma \rightarrow (ax_1+bx_2, cx_1+dx_2, z)\Gamma.$$

The group of isometries of the manifold $X$  can be seen as the
circle. Let $\theta$ be a Borel map defined from the quotient torus
$\TT^2$ to $\RR/\ZZ$: $\theta(x_1,x_2,z)=\theta(x_1,x_2)$. Then:
\begin{eqnarray*} & &\mu\{x: \
\frac{1}{\sqrt{n}}
\sum_{k=0}^{n-1}\varphi(A^k(x_1,x_2),\theta(x_1,x_2)+z) <t \}\\
& &\ \ \ \ \ \ \ \ \ \ \longrightarrow \int_{\RR/\ZZ}\frac{1}{\sqrt{2\pi} \sigma_y}
\int_{-\infty}^t \exp(-{1\over 2} \frac{s^2}{\sigma_y^2})ds\
d\mu_{\theta}(y),
\end{eqnarray*}
for every {H\"older} function $\varphi$,
where $\sigma_y^2$ is the asymptotic variance associated to the
function $\varphi_y(x_1,x_2,z)=\varphi(x_1,x_2,z+y)$.

\vskip 3mm{\bf Diagonal flows on  compact quotients of
$SL(d,\RR)$} ({\it cf.} \cite{LBP05})

Let $G$ be the group $SL(d,\RR)$, let $\Gamma$ be a cocompact
lattice of $G$, and let $\mu$ be the probability on $G/\Gamma$
deduced from the Haar measure. Let $g_0$ be a diagonal matrix in $G$
different from the identity. It defines a transformation $T$ on
$G/\Gamma$: $x =g\Gamma \longmapsto Tx = g_0 g\Gamma$.

\vskip 3mm Let $\varphi$ be a centered $C^\infty$ function from
$G/\Gamma$ to $\RR$ and let $\theta$ be a Borel map from $G/\Gamma$
to $SO(d,\RR)$. We have
\begin{eqnarray*} \mu\{x: \
\frac{1}{\sqrt{n}} \sum_{k=0}^{n-1}\varphi(\theta(x)T^kx) <t \}
\longrightarrow \int_{SO(d,\RR)}\frac{1}{\sqrt{2\pi} \sigma_y}
\int_{-\infty}^t \exp(-{1\over 2} \frac{s^2}{\sigma_y^2})ds\
d\mu_{\theta}(y),
\end{eqnarray*}
where $\sigma_y^2$ is the asymptotic variance associated to the
function $\varphi_y(x)=\varphi(y.x)$.

\vskip 3mm
\section{Generalization to multiple decorrelation \label{multiMix}}

\vskip 3mm \subsection{Multiple decorrelation and gaussian laws}

Let $(X,T,\mu)$ be a dynamical system defined on a manifold $X$. Let
$d$ be a Riemannian distance on $X$. The {H\"older} norm is defined as
in (\ref{etaNorm}). The expectation $\EE$ is the integral with
respect to $\mu$. In some proofs, we will denote by the same letter
$C$ a constant which may vary in the proof. Now we introduce the
following multiple decorrelation property:

\begin{ppty} \label{propH} There exist $C>0$ and $\delta\in]0,1[$
such that, for all integers $m$ and $m'$, all {H\"older} continuous
functions $(\varphi_i)_{i=1}^{m+m'}$ defined on $X$, all integers $0
\leq \ell_1 \leq \ldots\leq \ell_m \leq k_1 \leq \ldots\leq k_{m'},
N>0$,
$${\rm Cov}(\displaystyle\prod_{i=1}^{m}  T^{\ell_i }\,\varphi_i,
\displaystyle\prod_{j=1}^{m'} T^{k_j+N} \varphi_j )\leq
C(\displaystyle\prod_{i=1}^{m+m'}\| \varphi_i\|_{\infty}+
\displaystyle\mathop{\sum}_{j=1}^{m+m'}
[\varphi_j]_{\eta}\displaystyle\mathop{\prod}_{i\neq j}
\|\varphi_i\|_{\infty})\delta^N.$$ \end{ppty}

This property has many interesting consequences. Following C. Jan's
method \cite{Ja00}, we will show that the normalized ergodic sums
behave in some sense like gaussian variables even for some non
invariant measures on $X$ absolutely continuous with respect to
$\mu$ and which can vary with the time.

\vskip 3mm Let us first state a very simple consequence of this
property. There exist $C>0$, $\zeta \in ]0,1[$, such that, for every
centered {H\"older} continuous functions $\varphi$ and $\psi$ defined on
$X$ with zero average with respect to $\mu$, we have:
\begin{eqnarray}
|\langle \varphi,T^n\psi\rangle| \leq  C \|\varphi\|_\eta \,
\|\psi\|_2 \, \zeta^{n}, \ |\langle \varphi,T^n\psi\rangle| \leq  C
\|\varphi\|_2 \, \|\psi\|_\eta \, \zeta^n. \label{maj-dec}
\end{eqnarray}

These inequalities are consequences of  property \ref{propH}, proved by using the Cauchy-Schwarz inequality and
distinguishing two cases $\|\varphi\|_2\leq \|\varphi\|_\eta
\delta^{n/2}$ and $\|\varphi\|_2\geq \|\varphi\|_\eta \delta^{n/2}$.
In particular the asymptotic variance of the normalized ergodic sums
of a {H\"older} continuous function is well defined and given by:
$$\sigma(\varphi)^2=\mu(\varphi^2)+2\sum_{k=1}^{\infty}\langle\varphi,T^k\varphi\rangle.$$

\begin{thm} \label{tcl-densit} Let $(X,T,\mu)$ be a dynamical system defined
on a manifold $X$ satisfying Property \ref{propH}. Let $(\rho_n)$ be
a sequence of density functions with respect to $\mu$ with norms
$\|\rho_n\|_\eta$ bounded by $Cn^L$, for some constants $\eta
>0, C, L$. Then, for every sequence $(\varphi_n)$ of centered {H\"older}
continuous functions with {H\"older} norm uniformly bounded, we have,
for every real number $t$,
$$\mu(\rho_n \, \exp (\frac{it}{n^{1/2}}
\displaystyle\mathop{\sum}^{n-1}_{\ell=0}
T^{\ell}\varphi_n))-\exp(-{1\over 2} \sigma(\varphi_n)^2t^2)
\rightarrow 0.$$
\end{thm} \Proof \  \rm Let $\alpha\in]0,1/2[$. The difference
$\frac{1}{\sqrt{n}} S_n\varphi_n(x) - \frac{1}{\sqrt{n}}
S_n\varphi_n(x)\circ T^{n^\alpha}$ tends uniformly to 0. If
$\|\rho_n\|_\eta \leq Cn^L$, Property \ref{propH} gives:
\begin{eqnarray*}
& & |\mu(\rho_n \, \exp (\frac{it}{n^{1/2}}
\displaystyle\mathop{\sum}^{n-1}_{\ell=0}
T^{\ell+n^\alpha}\varphi_n))- \mu(\rho_n)\, \mu(\exp
\frac{it}{n^{1/2}} \displaystyle\mathop{\sum}^{n-1}_{\ell=0}
T^{\ell+n^\alpha}\varphi_n)|\\
& &\ \ \ \ \ \ \ \ \ \  \leq C(\|\rho_n\|_\infty+n
[\varphi_n]_{\eta}
\|\rho_n\|_{\infty}+[\rho_n]_{\eta})\delta^{n^\alpha}\leq
Cn^{L+1}\delta^{n^\alpha}.
\end{eqnarray*}

Thus we only have to study $\mu(\exp (\frac{it}{n^{1/2}}
\displaystyle\mathop{\sum}^{n-1}_{\ell=0}
T^{\ell+n^\alpha}\varphi_n)) = \mu(\exp \frac{it}{n^{1/2}}
\displaystyle\mathop{\sum}^{n-1}_{\ell=0} T^{\ell}\varphi_n)$.

Consider $(\Omega,\PP)$ a probability space containing $(X, \mu)$
and  a sequence $(X_{k,n})$ defined on $\Omega$ of centered
independent bounded random variables of variance
$\sigma^2(\varphi_n)$, independent from the variables $\varphi_n$,
with distribution
$1/2(\delta_{-\sigma(\varphi_n)}+\delta_{\sigma(\varphi_n)})$. Since
the sequence $(\|\varphi_n\|_\eta)_n$ is bounded, the sequence
$(\sigma(\varphi_n))_n$ is also bounded and one easily check that
$$ \EE( \exp (\frac{it}{n^{1/2}}
\displaystyle\mathop{\sum}^{n-1}_{\ell=0} X_{\ell,n}))-\exp(-{1\over
2}\sigma(\varphi_n)^2t^2)
=\cos^n({{1}\over{n^{1/2}}}\sigma(\varphi_n)t)-\exp(-{1\over 2}
\sigma(\varphi_n)^2t^2)$$
tends toward 0.

We claim that the difference between the characteristic functions of
$\frac{1}{n^{1/2}} \displaystyle\mathop{\sum}_{k=0}^{n-1}
T^k\varphi_n$ and $\frac{1}{n^{1/2}} {\sum}_{k=0}^{n-1} X_{k,n}$
tends to 0. To show it, let us write
$$B_{\ell,n} = \exp( \frac{it}{n^{1/2}} T^\ell\varphi_n),
\ \ \ C_{\ell,n} = \exp (\frac{it}{n^{1/2}} X_{\ell,n}).$$ One has:
\begin{eqnarray}
\exp \frac{it}{n^{1/2}} \displaystyle\mathop{\sum}^{n-1}_{\ell=0}
T^\ell\varphi_n- \exp \frac{it}{n^{1/2}}
\displaystyle\mathop{\sum}^{n-1}_{\ell=0}X_{\ell,n}
=\displaystyle\prod^{n-1}_{\ell=0}B_{\ell,n}-\displaystyle\prod^{n-1}_{\ell=0}
C_{\ell,n}, \label{BnCn}
\end{eqnarray}
 and
$$
\displaystyle\prod^{n-1}_{\ell=0} B_{\ell,n} -
\displaystyle\prod^{n-1}_{\ell=0} C_{\ell,n}=
\displaystyle\mathop{\sum}^{n-1}_{\ell=0}
(\displaystyle\prod^{\ell-1}_{k=0} C_{k,n})
\underbrace{(B_{\ell,n}-C_{\ell,n})
(\displaystyle\prod^{n-1}_{k=\ell+1} B_{k,n})}_{\Delta_{\ell}},
$$
where the products with an empty set of index are conventionally
taken to be 1.

The variables $\Delta_{\ell}=(B_{\ell,n}-
C_{\ell,n})\displaystyle\prod^{n^{1+\alpha}-1}_{k=\ell+1}
B_{\ell,n}$ and $\displaystyle\prod^{\ell-1}_{k=0} C_{\ell,n}$ are
independent. We will show that most of the $n$ terms
$|\EE(\Delta_{\ell})|$ are bounded by some constant times
$n^{-3/2}\ln n$. This will imply the result.

\vskip 3mm Consider a sequence  $\chi(m)$ that will be fixed later.
When $\ell+3 \chi(n)+1 < n$, we split the product $\Delta_{\ell}$ in
blocks:
\begin{eqnarray*}
\Delta_{\ell}&=&\underbrace{(B_{\ell,n}-C_{\ell,n})\displaystyle\prod^{\ell+\chi(n)}_{k=\ell+1}
B_{\ell,n}}_{{\cal A}}
\underbrace{\displaystyle\prod^{\ell+2\chi(n)}_{k=\ell+\chi(n)+1}B_{k,n}
}_{{\cal B}}
\underbrace{\displaystyle\prod^{\ell+3\chi(n)}_{k=\ell+2\chi(n)+1}
B_{k,n}}_{{\cal C}}
\underbrace{\displaystyle\prod^{n-1}_{k=\ell+3\chi(n)+1}
B_{k,n}}_{{\cal D}}.
\end{eqnarray*}
Let us write:
\begin{eqnarray}
\mathbb E(\Delta_{\ell})=\mathbb E({\cal ABCD})= \mathbb E( {\cal
A}({\cal B}-1) ({\cal C}-1){\cal D})+\mathbb E( {\cal ABD})+\mathbb
E( {\cal ACD}) -\mathbb E({\cal AD}).\label{decomp}
\end{eqnarray}
The mean value theorem shows that ${\cal A}$ is bounded by
$$C t n^{-1/2}(\|\varphi_n\|_{\infty}+\|X_{\ell,n}\|_{\infty})$$
for some constant $C$, and $({\cal B}-1), ({\cal C}-1)$ are both
bounded by
$$\frac{2t}{n^{1/2}}
\displaystyle\mathop{\sum}^{\ell+2\chi(n)}_{k=\ell+\chi(n)+1}
\|\varphi_n\|_{\infty}.$$

Thus $\mathbb E( {\cal A}({\cal B}-1) ({\cal C}-1) {\cal D}) \leq
C\|\varphi_n\|_{\infty}^2(\|\varphi_n\|_{\infty}+\|X_{\ell,n}\|_{\infty})
\frac{t^3}{n^{3/2}} \chi(n)^2$.

The sequences $\|\varphi_n\|_{\infty}$,  $\|\varphi_n\|_{\eta}$,
$\|X_{\ell,n}\|_{\infty}$ being  bounded, we will not retain the
dependence on these quantities in our computation, nor the
dependence on $t$. For example we just write
\begin{eqnarray}
\mathbb E( {\cal A}({\cal B}-1) ({\cal C}-1) {\cal D}) \leq  C
\frac{1}{n^{3/2}} \chi(n)^2.\label{maj1}
\end{eqnarray}

\vskip 3mm The three other terms in (\ref{decomp}) can be treated in
the following way.

Consider for example : $\mathbb E( {\cal ABD})=Cov ( {\cal
AB,D})+\mathbb E ( {\cal AB}) \mathbb E({\cal D})$. Property
\ref{propH} gives
\begin{eqnarray}
Cov ( {\cal AB, D})\leq C\  \zeta^{\chi(n)}.\label{maj3}
\end{eqnarray}
Let us now study
\begin{eqnarray*} {\cal
AB}&=&(B_{\ell,n}-C_{\ell,n})
\displaystyle\prod^{\ell+2\chi(n)}_{k=\ell+1}
B_{k,n}\\
&=& \left(\exp\frac{it  T^\ell\varphi_n}{n^{1/2}}-\exp \frac{it
X_{\ell,n}}{n^{1/2}}\right) \exp \left(itn^{-1/2}
\displaystyle\mathop{\sum}^{\ell+2\chi(n)}_{k=\ell+1}  T^k\varphi_n\right)\\
\end{eqnarray*}
The expansion of the two terms of this product at order 1 and order
2 yields
\begin{eqnarray*}
{\cal AB} &=&\frac{it}{n^{1/2}}
(T^\ell\varphi_n-X_{\ell,n})-\frac{t^2}{2n}
(T^\ell\varphi_n^2 -X_{\ell,n}^2)\\
&{ }& \ \ -\frac{t^2}{n}
\displaystyle\mathop{\sum}^{\ell+2\chi(n)}_{k=\ell+1} T^k\varphi_n \
T^\ell\varphi_n { }\ \ +\frac{t^2}{n}
X_{\ell,n}\displaystyle\mathop{\sum}^{\ell+2\chi(n)}_{k=\ell+1}
T^k\varphi_n { }  +D,
\end{eqnarray*}
where $D$ satisfies
\begin{eqnarray*}|D| \leq
C[\frac{t^2}{n^{3/2}}(\chi(n)^2+\frac{t^3}{n^{3/2}}
+\frac{t^4}{n^{2}}\chi(n)]\leq C {{1}\over{n^{3/2}}}\chi(n)^2.
\end{eqnarray*}
By taking the expectation, one obtains:
\begin{eqnarray}
&&|\EE( {\cal AB})| \leq  \nonumber \\
&&\frac{t^2}{2n}\left( \EE (X_{\ell,n}^2) - \left(\mathbb
E(T^{\ell}\varphi_n^2)+2
\displaystyle\mathop{\sum}^{\ell+2\chi(n)}_{k=\ell+1}\mathbb
E(T^k\varphi_n \ T^{\ell}\varphi_n)\right)\right) + C
{{\chi(n)^2}\over{n^{3/2}}}. \label{MajAB}
\end{eqnarray}
By definition of the $X_{\ell,n}$, we have
\begin{eqnarray*}
\EE(X_{\ell,n}^2)&=&\sigma^2(\varphi_n) = E(T^{\ell}\varphi_n^2)+2
\displaystyle\mathop{\sum}^{\infty}_{k=\ell+1}\mathbb
E(T^k\varphi_n \ T^{\ell}\varphi_n) \\
&=& E(T^{\ell}\varphi_n^2)+2
\displaystyle\mathop{\sum}^{\ell+2\chi(n)}_{k=\ell+1}\mathbb
E(T^k\varphi_n \ T^{\ell}\varphi_n) \ +2
\displaystyle\mathop{\sum}^{\infty}_{k=\ell+2\chi(n)+1}\mathbb
E(T^k\varphi_n \ T^{\ell}\varphi_n).
\end{eqnarray*}
Replacing $\EE(X_{\ell,n}^2)$ by this expression in (\ref{MajAB}) we
obtain
$$|\EE( {\cal AB})|\leq C(\frac{t^2}{2n}\displaystyle\mathop{\sum}^{\infty}_{k=\ell+2\chi(n)+1}\mathbb
E(T^k\varphi_n \ T^{\ell}\varphi_n) +\chi(n)^2n^{-3/2}),
$$
and, because of the decay of correlations (\ref{maj-dec}):
\begin{eqnarray}
|\EE( {\cal AB})|\leq C(
{{\zeta^{\chi(n)}}\over{n}}+\chi(n)^2n^{-3/2}).\label{maj2}
\end{eqnarray}
Since $|\mathbb E({\cal D})|\leq 1$, the inequalities (\ref{maj1}),
(\ref{maj3}) and (\ref{maj2}) give (the constant $C$ may change):
$${\mathbb E}(\Delta_{\ell})\leq
C({{\zeta^{\chi(n)}}\over{n}}+\chi(n)^2n^{-3/2}+\zeta^{\chi(n)}).$$

\vskip 3mm Now we can bound (\ref{BnCn}):
\begin{eqnarray*}
&&|\mathbb E(\displaystyle\prod^{n-1}_{0} B_{\ell,n} -
\displaystyle\prod^{n-1}_{0} C_{\ell,n})|=
|\displaystyle\mathop{\sum}^{n-1}_{\ell=0} \mathbb
E(\displaystyle\prod^{\ell-1}_{k=0} C_{k,n}) \mathbb
E(\Delta_{\ell})| \leq \displaystyle\mathop{\sum}^{n-1}_{\ell=0}
 |\mathbb E(\Delta_{\ell})|\\
&\leq&\displaystyle\mathop{\sum}^{n-1-3\chi(n)-1}_{\ell=0}C({{\zeta^{\chi(n)}}
\over{n}}+\chi(n)^2n^{-3/2}+\zeta^{\chi(n)})+\displaystyle\mathop{\sum}^{n-1}
_{\ell=n-1-3\chi(n)}
{\mathbb E}(\Delta_{\ell}).
\end{eqnarray*}

The mean value theorem implies $\mathbb E(\Delta_{\ell})\leq
Cn^{-1/2}$.

If we take $\chi(n)=D\ln n$ with $D>-\ln(\zeta)^{-1}$, then
(\ref{BnCn}) tends to 0, in view of
\begin{eqnarray*}
|\mathbb E(\displaystyle\prod^{n-1}_{0} B_{\ell,n} -
\displaystyle\prod^{n-1}_{0} C_{\ell,n})|&\leq&C(\zeta^{D\ln
n}+\frac{D^2\ln^2(n)}{n^{1/2}}+n\zeta^{D\ln(n)}+\frac{3D\ln(n)}{n^{1/2}}).
\end{eqnarray*}
\fdm

\vskip 3mm
\subsection{Application}

Let $(X,T,\mu)$ be a dynamical system defined on a manifold $X$.
Suppose that Property \ref{propH} holds. Let $\theta : x \rightarrow
\theta_x$ be a map from $X$ into the set of the {H\"older} continuous
maps from $X$ to $X$. We will use also the notation $\theta(x)$
instead of $\theta_x$. We suppose that $\theta$ satisfies for some
$\eta
> 0$
\begin{eqnarray}
\sup_{z\in X}d(\theta_x(z),\theta_y(z))\leq Cd(x,y)^\eta.
\label{regTheta}
\end{eqnarray}

Let $\varphi$ be a {H\"older} continuous function from $X$ to $\RR$. We
want to study the behavior of
$$\sum_{k=0}^{n-1}\varphi(\theta_x(T^k x)).$$

We can assume that the exponent of regularity in (\ref{regTheta})
and the exponent for $\varphi$ are the same. Since the maps
$\theta_x $ do not necessarily preserve the measure $\mu$, one has
to center these modified ergodic sums. We write
$${\overline \varphi_{\theta_x }} := \EE(\varphi(\theta_x (.))) =
\int_X \varphi(\theta_x (y)) d\mu(y),$$
$$\varphi_{\theta_x } := \varphi(\theta_x \cdot)-{\overline
\varphi_{\theta_x }}, \ \ \sigma^2_{\theta_x} :=
\sigma^2(\varphi_{\theta_x }).$$

For $\mu$-almost every $x$ (the generic points of the dynamical
system),
$${{1}\over{n}}\sum_{k=0}^{n-1}\varphi(\theta_x (T^kx)) \rightarrow
{\overline \varphi_{\theta_x }}.$$

We claim that both functions ${\overline \varphi_{\theta_x }}$ and
$\sigma^2_{\theta_x }$ are {H\"older} continuous. Let us show it for
$\sigma^2_{\theta_x }$. On one hand, in view of (\ref{maj-dec}), the
absolute value of
\begin{eqnarray*}
\sigma^2(\varphi_{\theta_x })-\sigma^2(\varphi_{\theta_y }) =
\mu(\varphi_{\theta_x }^2)-\mu(\varphi_{\theta_y
}^2)+2\sum_1^{\infty} (\mu(\varphi_{\theta_x }.T^k\varphi_{\theta_x
})-\mu(\varphi_{\theta_y }.T^k\varphi_{\theta_y }))
\end{eqnarray*}
is bounded by
\begin{eqnarray*}
&&2(\|\varphi\|_\infty\|\varphi_{\theta_x }-\varphi_{\theta_y }\|_2
+\sum_1^{\infty}|\mu(\varphi_{\theta_x }(T^k\varphi_{\theta_x
}-T^k\varphi_{\theta_y })|
+\sum_1^{\infty}| \mu((\varphi_{\theta_x }-\varphi_{\theta_y })T^k\varphi_{\theta_y })|)\\
&&\leq 2\|\varphi\|_\infty\|\varphi_{\theta_x }-\varphi_{\theta_y
}\|_2 +4C\sum_1^{\infty}C(\|\varphi_{\theta_x }\|_\eta
+\|\varphi_{\theta_y }\|_\eta)\|\varphi_{\theta_x }-\varphi_{\theta_y }\|_2 \, \zeta^k\\
&&\leq C(\|\varphi_{\theta_x }\|_\eta+\|\varphi_{\theta_y
}\|_\eta)\|\varphi_{\theta_x }-\varphi_{\theta_y }\|_2,
\end{eqnarray*}
and, on the other hand,
$$\|\varphi_{\theta_x }-\varphi_{\theta_y }\|_2\leq \|\varphi_{\theta_x }-\varphi_{\theta_y }\|_\infty\leq
Cd(x,y)^{\eta^2}.$$

Now suppose that there exists $({\cal P}_n)$ a sequence of
partitions with the diameter of the elements of ${\cal P}_n$ smaller
than $n^{-2/\eta^2}$ such that, for positive constants $C,L$, for
every $n \in \NN$, every $P$ in ${\cal P}_n$, there exists a density
function $\rho_{n,P}$ such that
\begin{eqnarray}
\|\rho_{n,P}\|_\eta \leq Cn^L, \ \|\rho_{n,P}- \mu(P)^{-1}{\bf
1}_P\|_1 \leq {1 \over n}. \label{density-part}
\end{eqnarray}

Let us fix such a sequence $({\cal P}_n)$. Cutting out the space $X$
according to the partition ${\cal P}_n$ we get
$$\sum_{k=0}^{n-1}\varphi_{\theta_x }(T^kx)=\sum_{P\in{\cal
P}_n}\sum_{k=0}^{n-1}{\bf 1}_P(x)\varphi_{\theta_x }(T^kx).$$

For each element $P$ of ${\cal P}_n$ we choose in $P$ a point
$x_{n,P}$. If $x\in P$, then the distance $d(\theta_x (T^k
x),\theta_{x_{n,P}}(T^k x))$ is bounded by $Cd(x,x_{n,P})^\eta\leq
Cn^{-2/\eta}$. We thus have
$$|\sum_{P\in{\cal P}_n}\sum_{k=0}^{n-1}{\bf 1}_P(x)\varphi_{\theta_x
}(T^kx)-\sum_{P\in{\cal P}_n}\sum_{k=0}^{n-1}{\bf
1}_P(x)\varphi_{\theta_{x_{n,P}}}(T^kx)|\leq C n n^{-2}.$$

Let us now study the characteristic function of $\sum_{P\in{\cal
P}_n}\sum_{k=0}^{n-1}{\bf 1}_P(x)\varphi_{\theta(x_{n,P})}(T^kx)$:
\begin{eqnarray*}
& &\EE(\exp(i {t\over \sqrt{n}} \sum_{P\in{\cal
P}_n}\sum_{k=0}^{n-1}{\bf 1}_P(x)\, \varphi_{\theta(x_{n,P})}(T^kx)))\\
& &\ \ \ \ \ \ =\EE(\sum_{P\in{\cal P}_n}{\bf
1}_P(x)\exp(i {t\over \sqrt{n}} \sum_{k=0}^{n-1}\varphi_{\theta(x_{n,P})}(T^kx)))\\
& &\ \ \ \ \ \ =\sum_{P\in{\cal P}_n}\mu(P)\EE(\mu(P)^{-1}{\bf
1}_P(x) \, \exp(i {t\over \sqrt{n}}
\sum_{k=0}^{n-1}\varphi_{\theta(x_{n,P})}(T^kx))).
\end{eqnarray*}

By (\ref{density-part}) we have
$$
|\EE(\frac{{\bf 1}_P(x)}{\mu(P)} \, \exp(i {t\over \sqrt{n}}
\sum_{k=0}^{n-1}\varphi_{\theta(x_{n,P})}
(T^kx)))-\EE(\rho_{n,P}\exp(i {t\over \sqrt{n}} \sum_{k=0}^{n-1}
\varphi_{\theta(x_{n,P})} (T^kx)))|\leq 1/n.$$

We know from Theorem \ref{tcl-densit} that
$$\EE(\rho_{n,P} \, \exp(i{t\over \sqrt{n}}
\sum_{k=0}^{n-1} \varphi_{\theta(x_{n,P})} (T^kx))) - \exp(-{1\over
2} \, \sigma^2_{ \theta(x_{n,P})}t^2) \rightarrow 0.$$

The sum $\sum_{P\in{\cal P}_n}\mu(P)\ \exp(-{1\over 2} \sigma^2_{
\theta(x_{n,P})}t^2)$ is a Riemann sum of the {H\"older} continuous
function $\sigma^2_{\theta_x }$. It converges to
$$\int_X \exp(- {1\over 2} \sigma^2_{\theta_x }t^2)\ d\mu(x).$$

Now for the assumption (\ref{density-part}), it is known that every
smooth manifold has a triangulation (see \cite{Ca35,Wh57}).
The existence of a sequence of partitions ${\cal P}_n$ and of
functions $\rho_{n,P}$ satisfying the previous conditions are thus
satisfied when $X$ is a smooth manifold. So, we have proved the
following theorem.

\begin{thm}
Let $(X,T,\mu)$ be a dynamical system defined on a smooth manifold
$X$ for which Property \ref{propH} holds. Let $\theta$ be a map from
$X$ to the set of the {H\"older} continuous maps from  $X$ to $X$ such
that (\ref{regTheta}) is satisfied. Then for every centered {H\"older}
continuous function $\varphi$ on $X$, we have
$$\EE [\exp(i {t \over \sqrt{n}}\sum_{k=0}^{n-1}(\varphi(\theta_x T^kx)
-{\overline \varphi_{\theta_x }}))] \rightarrow \int_X \exp(-
{1\over 2} \sigma^2_{\theta_x }t^2)\ d\mu(x).$$
\end{thm}

Now let us give some examples satisfying the hypothesis. The
examples given at the end of Section \ref{groupAct} can be treated
by this method since Property \ref{propH} is satisfied for ergodic
automorphisms of the torus, ergodic automorphisms on nilmanifolds,
diagonal flows on  compact quotients of $SL(d,\RR)$. By the method
used in this Section \ref{multiMix} we obtain the same result for
{H\"older} continuous functions under the assumption that the map
$\theta$ is {H\"older} continuous. If we identify in these cases the
element $\theta_x $ of $K$ with the translation by $\theta(x) $ we
obtain:
$$\EE(\exp(i {t \over \sqrt{n}}\sum_{k=0}^{n-1}\varphi(\theta(x)(T^kx)))
\rightarrow \int_X \exp(-{1\over 2} \sigma^2_{\theta(x)}t^2)\
d\mu(x)$$
and, by definition of the image measure of $\mu$ by
$\theta$ on $K$:
$$\int_X \exp(- {1\over 2} \sigma^2_{\theta(x)}t^2)\ d\mu(x)=\int_K
\exp(- {1\over 2} \sigma^2_{k}t^2)\ d\mu_{\theta}(k).$$

By linearity of the Fourier transform, we just have another
formulation of the convergence stated in Proposition \ref{groupes}.
The result is the same. But this second method allows us to consider
transformations $\theta(x)$ that are not given by translation by
elements of a compact group acting on $X$. In the case of the
automorphisms of the torus, one can for example take for $\theta_x$
a regular family of diffeomorphisms of $\TT^d$. This method can be
also adapted to expanding dynamical systems for which the
Perron-Frobenius operator has good spectral properties.

\subsection{Convergence of the variance \label{subsect-conv-var}}

\begin{prop} Let $(X,T,\mu)$ be a dynamical system defined on a manifold
$X$ for which Property \ref{propH} holds. Let $\theta$ be a map from
$X$ to the set of the {H\"older} continuous maps from  $X$ to $X$ such
that (\ref{regTheta}) is satisfied, and $\varphi$ a centered {H\"older}
continuous function  on $X$. Then
$$\EE\left({{1}\over{n}}\left(\sum_{k=1}^n\varphi(\theta(\cdot)T^k\cdot)
\right)^2\right)\rightarrow \int_X\sigma^2_{\theta_x }\ d\mu(x).$$
\end{prop}
\Proof \ \ Consider a sequence of partitions $({\cal P}_n)$ with the
same properties as in the previous subsection. We have:
\begin{eqnarray*}
\EE\left({{1}\over{n}}\left(\sum_{k=1}^n\varphi(\theta(\cdot)T^k\cdot)
\right)^2\right)&=&\EE\left(\sum_{P\in{\cal P}_n}{\bf
1}_P(\cdot){{1}\over{n}}\left(\sum_{k=0}^{n-1}\varphi_{\theta(\cdot)} (T^k\cdot)\right)^2\right)\\
&=&\sum_{P\in{\cal P}_n}\mu(P) \, \EE\left({{{\bf
1}_P(\cdot)}\over{\mu(P)}} {{1}\over{n}}
\left(\sum_{k=0}^{n-1}\varphi_{\theta(\cdot)}(T^k\cdot)\right)^2\right).
\end{eqnarray*}
For every $P\in {\cal P}_n$, let $\rho_{n,P}$ be a function with
norm $\|\rho_{n,P}\|_\eta$ less than $Cn^L$ such that $\|\rho_{n,P}-
\mu(P)^{-1}{\bf 1}_P\|_1\leq n^{-2}$. We have:
$$|\EE\left({{1}\over{n}} \left(\sum_{k=1}^n\varphi(\theta(\cdot)T^k\cdot)\right)^2\right) -
\sum_{P\in{\cal P}_n} \mu(P)\, \EE\left(\rho_{n,P}
{{1}\over{n}}\left(\sum_{k=0}^{n-1}\varphi_{\theta(x_{n,P})}
(T^k\cdot)\right)^2\right)|\leq 1/n.$$

Let us study
$\EE\left(\rho_{n,P}{{1}\over{n}}\left(\sum_{k=0}^{n-1}\varphi_{\theta(x_{n,P})}(T^k\cdot)\right)^2\right)$.
We fix $\alpha\in(0,1/2)$.
\begin{eqnarray*}
&&\left|\EE(\rho_{n,P}{{1}\over{n}}(\sum_{k=0}^{n-1}\varphi_{\theta(x_{n,P})}(T^k\cdot))^2)
-\EE(\rho_{n,P}{{1}\over{n}} (\sum_{k=n^\alpha}^{n+n^\alpha-1}\varphi_{\theta(x_{n,P})}(T^k\cdot))^2) \right|\\
&&=\EE(\rho_{n,P}{{1}\over{n}}(\sum_{k=0}^{n^\alpha-1}\varphi_{\theta(x_{n,P})}
(T^k\cdot)\\
& & \ \ \ \ \ \ +\sum_{k=n^\alpha}^{n-1}\varphi_{\theta(x_{n,P})}(T^k\cdot))^2)
-(\sum_{k=n^\alpha}^{n-1}\varphi_{\theta(x_{n,P})}
(T^k\cdot)+\sum_{k=n}^{n+n^\alpha-1}\varphi_{\theta(x_{n,P})}(T^k\cdot))^2)\\
&&\  ={{1}\over{n}}\EE(\rho_{n,P} [(\sum_{k=0}^{n^\alpha-1}
\varphi_{\theta(x_{n,P})}(T^k\cdot))^2+(\sum_{k=n}^{n+n^\alpha-1}\varphi_{\theta(x_{n,P})}(T^k\cdot))^2)]\\
& & \ \ \ \ \ \ +2 \EE(\rho_{n,P} (\sum_{k=n^\alpha}^{n-1}
\varphi_{\theta(x_{n,P})}(T^k\cdot))(\sum_{k=0}^{n^\alpha-1}\varphi_{\theta(x_{n,P})}
(T^k\cdot)-\sum_{k=n}^{n+n^\alpha-1}\varphi_{\theta(x_{n,P})}(T^k\cdot)))\\
& &\ \ \ \ \leq n^{2\alpha-1}\|\varphi\|^2_\infty+{{4}\over{n}}
[\EE(\rho_{n,P}(\sum_{k=n^\alpha}^{n-1}
\varphi_{\theta(x_{n,P})}(T^k\cdot))^2)]^{1/2}
n^\alpha\|\varphi\|_\infty
\end{eqnarray*}

We have
\begin{eqnarray*}
&& \EE(\rho_{n,P}(\sum_{k=n^\alpha}^{n-1}
\varphi_{\theta(x_{n,P})}(T^k\cdot))^2) \\
&& \leq \EE(\rho_{n,P}) \, \EE((\sum_{k=n^\alpha}^{n-1}
\varphi_{\theta(x_{n,P})}(T^k\cdot))^2) +|Cov(\rho_{n,P},
(\sum_{k=n^\alpha}^{n-1}\varphi_{\theta(x_{n,P})}(T^k\cdot))^2)|.
\end{eqnarray*}

Because of Property \ref{propH} the covariance tends to 0 faster
than $1/n$. On the other hand
$\EE((\sum_{k=n^\alpha}^{n-1}\varphi_{\theta(x_{n,P})}(T^k\cdot))^2)$
is bounded by $Cn$ for some constant $C>0$. Putting this together
one obtains:
\begin{eqnarray*}
&
&|\EE(\rho_{n,P}{{1}\over{n}}(\sum_{k=0}^{n-1}\varphi_{\theta(x_{n,P})}
(T^k\cdot))^2)-\EE(\rho_{n,P}{{1}\over{n}}(\sum_{k=n^\alpha}^{n+n^\alpha-1}
\varphi_{\theta(x_{n,P})}(T^k\cdot))^2)|\\
& & \ \ \ \ \ \ \ \ \ \ \ \leq
n^{2\alpha-1}\|\varphi\|^2_\infty+4Cn^{\alpha-1}\|\varphi\|_\infty
\|\varphi\|_\eta n^{1/2} \leq C\|\varphi\|_\infty\|\varphi\|_\eta
n^{\alpha-1/2}.
\end{eqnarray*}
Again because of Property \ref{propH}, we have
$$|\EE(\rho_{n,P}{{1}\over{n}}(\sum_{k=n^\alpha}^{n+n^\alpha-1}\varphi_{\theta(x_{n,P})}(T^k\cdot))^2)
-\EE({{1}\over{n}}(\sum_{k=n^\alpha}^{n+n^\alpha-1}\varphi_{\theta(x_{n,P})}(T^k\cdot))^2)|\leq
C\|\varphi\|_\infty\|\varphi\|_\eta n^{-1}$$

Thus we have shown that
$$ \left|\EE({{1}\over{n}}(\sum_{k=1}^n\varphi(\theta(\cdot)T^k\cdot))^2)
- \sum_{P\in{\cal P}_n} \mu(P) \,
\EE({{1}\over{n}}(\sum_{k=0}^{n-1}\varphi_{\theta(x_{n,P})}
(T^k\cdot))^2)\right|\leq C\|\varphi\|_\infty\|\varphi\|_\eta
n^{\alpha-1/2}.$$

Because of the mixing property of $T$ we have
$$\left|\EE({{1}\over{n}}(\sum_{k=0}^{n-1}\varphi_{\theta(x_{n,P})}(T^k\cdot))^2)
- \sigma^2_{\theta(x_{n,P})}\right|\leq
C\|\varphi\|_\infty\|\varphi\|_\eta n^{-1},$$
so that
$$\left|\EE({{1}\over{n}}(\sum_{k=1}^n\varphi(\theta(\cdot)T^k\cdot))^2)-\sum_{P\in{\cal
P}_n}\mu(P) \, \sigma^2_{\theta(x_{n,P})}\right|\leq
C\|\varphi\|_\infty\|\varphi\|_\eta n^{\alpha-1/2}.$$

Since the last sum is a Riemann sum for the function $x\rightarrow
\sigma^2_{\theta_x }$, this concludes the proof. \fdm

\vskip 3mm
\bibliographystyle{amsalpha}

\begin{thebibliography}{12345678}

\bibitem[1]{AiBe08}
Aistleitner (C.), Berkes (I.): On the central limit theorem for
$f(n_k x)$, (2008) to appear.

\bibitem[2]{Be76}
Berkes (I.): On the asymptotic behaviour of $Sf(n_{k}x)$. Main
theorems. Z. Wahrscheinlichkeitstheorie und Verw. Gebiete 34 (1976),
no. 4, p. 319-345.

\bibitem[3]{Bo82}
Bourbaki (N.): El\'ements de math\'ematiques: groupes et alg\`ebres de
Lie, Chap. 9, Groupes de Lie r\'eels compacts, Masson, Paris (1982).

\bibitem[4]{Ca35}
Cairns (S. S.): {Triangulation of the manifold of class one}, {Bull.
Amer. Math. Soc.}, {41} (1935), {8}, p. 549-552.

\bibitem[5]{CoLB02}
Conze (J.-P.), Le Borgne (S.): Le TCL pour une classe
d'automorphismes non hyperboliques de nilvari\'et\'es (2002) \hfill
\break http://perso.univ-rennes1.fr/stephane.leborgne/.

\bibitem[6]{Ea76}
Eagleson  (G. K.): {Some simple conditions for limit theorems to be
mixing}, {Teor. Verojatnost. i Primenen.}, {21} (1976) {3}, p.
653-660.

\bibitem[7]{FeFo47}
Ferrand (J.), Fortet (R.): Sur des suites arithm\'etiques
\'equir\'eparties, C. R. Acad. Sci. Paris 224 (1947), p. 516-518.

\bibitem[8]{Fo40}
Fortet (R.): Sur une suite \'egalement r\'epartie, Studia Math. 9
(1940), p. 54-70.

\bibitem[9]{Ga66}
Gaposhkin (V. F.): Lacunary series and independent functions.
(Russian) Uspehi Mat. Nauk 21 (1966) no. 6 (132), p. 3-82.

\bibitem[10]{Ja00}
Jan (C.): {Vitesse de convergence dans le {TCL} pour des cha\^\i nes
de {M}arkov et certains processus associ\'es \`a des syst\`emes
dynamiques}, {C. R. Acad. Sci. Paris S\'er. I Math.}, {331} (2000)
{5}, p. 395-398.

\bibitem[11]{Ka46}
Kac (M.): On the distribution of values of sums of the type $\sum
f(2\sp k t)$. Ann. of Math. (2) 47 (1946), p. 33-49.

\bibitem[12]{Ka49}
Kac (M.): Probability methods in some problems of analysis and
number theory. Bull. Amer. Math. Soc. 55 (1949), p. 641-665.

\bibitem[13]{Ka47}
Kac (M.): Distribution properties of certain gap sequences, Bull.
Amer. Math. Soc. Abstract 53, 7, p. 290.

\bibitem[14]{LBP05}
{Le Borgne (S.), P{\`e}ne (F.): {Vitesse dans le th\'eor\`eme limite
central pour certains syst\`emes dynamiques quasi-hyperboliques},
{Bull. Soc. Math. France}, {133} (2005) 3, p. 395-417.}

 \bibitem[15]{LB99}
Le Borgne (S.): {Limit theorems for non-hyperbolic automorphisms of
the torus}, {Israel J. Math.}, {109} (1999), p. 61-73.

\bibitem[16]{Le60}
Leonov  (V. P.): {On the central limit theorem for ergodic
endomorphisms of compact commutative groups}, {Dokl. Akad. Nauk
SSSR}, {135} (1960), p. 258-261.

\bibitem[17]{SaZy48}
Salem (R.), Zygmund (A.): On lacunary trigonometric series. II.
Proc. Nat. Acad. Sci. U. S. A. 34 (1948), p. 54-62.

\bibitem[18]{Ve86}
Veech (W.): Periodic points and invariant pseudomeasures for toral
endomorphisms. Ergodic Theory Dynam. Systems 6 (1986), no. 3, p.
449-473.

\bibitem[19]{Wh57}
Whitney (H.): {Geometric integration theory}, {Princeton University
Press} (1957).

\bibitem[20]{Zw07}
Zweim\"uller (R.): Mixing limit theorems for ergodic transformations,
J. Theor. Prob. (2007) 20, p. 1059-1071.

\end{thebibliography}

\end{document}